\documentclass[matann,numbook,envcountsame]{svjour}

\usepackage{amsmath,latexsym,amssymb,times,mathptm}

\def\sqr#1#2{{\vcenter{\hrule height.#2pt
        \hbox{\vrule width.#2pt height#1pt \kern#1pt
                \vrule width.#2pt}
        \hrule height.#2pt}}}
\def\square{\mathchoice\sqr64\sqr64\sqr{4}3\sqr{3}3}
\def\QED{\hfill$\square$}

\def\tratto{\mbox{\rule{2mm}{.2mm}$\;\!$}}

\begin{document}

\title{The structure of the core of ideals\thanks{The second and third
authors were partially supported by the NSF}}

\dedication{Dedicated to the memory of Professor Manfred Herrmann}

\titlerunning{The structure of the core of ideals}
\authorrunning{A. Corso et al.}

\author{Alberto Corso \and Claudia Polini \and Bernd Ulrich}


\institute{{\sc A. Corso} \\
Department of Mathematics, University of Kentucky, Lexington, KY
40506, USA \\ (e-mail$\colon$ corso@ms.uky.edu) \and
{\sc C. Polini} \\
Department of Mathematics, University of Oregon, Eugene, OR 97403,
USA \\ (e-mail$\colon$ polini@math.uoregon.edu) \and
{\sc B. Ulrich} \\
Department of Mathematics, Michigan State University, East
Lansing, MI 48824, USA \\ (e-mail$\colon$ ulrich@math.msu.edu)}

\date{Received: May 16, 2000 / Revised version: December 11, 2000}

\maketitle

\begin{abstract}
The core of an $R$-ideal $I$ is the intersection of all
reductions of $I$. This object was introduced by D. Rees and J. Sally and
later studied by C. Huneke and I. Swanson, who showed in particular
its connection to J. Lipman's notion of adjoint of an ideal.

Being an a priori infinite intersection of ideals, the core is
difficult to describe explicitly. We prove in a broad setting
that: ${\rm core}(I)$ is a finite intersection of minimal
reductions; ${\rm core}(I)$ is a finite intersection of general
minimal reductions; ${\rm core}(I)$ is the contraction to $R$ of a
`universal' ideal; ${\rm core}(I)$ behaves well under flat
extensions. The proofs are based on general multiplicity estimates
for certain modules.

\medskip

\noindent \subclassname{Primary 13A30, 13B21, 13H15; Secondary
13C40, 13H10.}
\end{abstract}

\section{Introduction}

To study properties of an ideal $I$ of a Noetherian local ring $R$
one often passes to a different ideal, either larger
or smaller, which carries most of the information about the
original ideal, but has better features.
The radical $\sqrt{I}$ of $I$, and the importance
the Nullstellensatz assigns to it, is the most notorious example.
The integral closure $\overline{I}$ or a reduction $J$ of $I$
are also very familiar instances.

We recall that the {\it radical} $\sqrt{I}$ of $I$ consists of all
solutions in $R$ of equations of the form $X^m-a=0$, with $a\in I$
and $m$ a non negative integer. The {\it integral closure}
$\overline{I}$ of $I$ consists instead of all solutions in $R$ of
equations of the form $X^n+b_1 X^{n-1} + b_2 X^{n-2} + \ldots +
b_{n-1} X + b_n = 0$, with $b_j \in I^j$ and $n$ a non negative
integer. We clearly have $I \subset \overline{I} \subset
\sqrt{I}$, with, in general, strict inclusions. Finally, a {\it
reduction} $J$ of $I$ is a subideal of $I$ such that $\overline{J}
= \overline{I}$. Equivalently, $J \subset I$ is a reduction of $I$
if $I^{r+1} = JI^r$ for some non negative integer $r$. Minimal
reductions are reductions which are minimal with respect to
containment. If the residue field of the ring $R$ is infinite,
then minimal reductions have the same number of generators, namely
the {\it analytic spread} $\ell = \ell(I)$ of $I$.

A more familiar description is the one of $\sqrt{I}$ as the intersection of
all prime ideals containing $I$ or, equivalently, as the intersection of all
minimal primes over $I$. It is well known that this intersection
is finite. Also, by work of D. Eisenbud, C. Huneke and W.V. Vasconcelos
\cite{EHV}, it is now easy to give an algorithmic approach to
$\sqrt{I}$ suitable for effective computer calculations.

On the other hand, reductions of an ideal are highly non unique.
Their intersection, dubbed {\it core} of the ideal $I$, comes from a
more recent vintage. It was studied for
the first time by D. Rees and J. Sally \cite{RS} and later by C. Huneke
and I. Swanson \cite{HS}, who also showed a connection with work of
J. Lipman on the adjoint of an ideal \cite{L}.
Being the intersection of an {\it a priori} infinite number of ideals,
this object is difficult to describe in terms of explicit data attached
to the ideal. It is known though that $\sqrt{{\rm core}(I)} = \sqrt{I} $, see
\cite{Vas} for instance.

The core of an ideal appears naturally in the context of the Brian\c{c}on--Skoda
theorem \cite{LS}. In one of its simplest formulations, this theorem says that if $R$
is a regular local ring of dimension $d$ and $I$ is an ideal then
$\overline{I^d} \subset J$ for every reduction $J$ of $I$, or equivalently,
$\overline{I^d} \subset {\rm core}(I)$.

The issues we address in this paper and to which we give fairly
general affirmative answers are: Is the core a finite intersection
of minimal reductions of $I$? Is the core a finite intersection of
general minimal reductions of $I$? Is the core the contraction to
$R$ of a `universal $\ell$-generated ideal'? Does the core behave
well under flat extensions? The last question has already been
raised by C. Huneke and I. Swanson in \cite{HS}.

Our results are based on general multiplicity estimates for
certain modules $($Lemmas~\ref{4.1},~\ref{4.2} and~\ref{4.4}$)$,
and their proofs use techniques coming from the theory of residual
intersections. We are required to introduce and base our
constructions on the notions of generic, universal and general
ideals. To be more specific, let $(R, {\mathfrak m})$ be a local
Cohen--Macaulay ring with infinite residue field and $I=(f_1,
\ldots, f_n)$ an $R$-ideal of height $g$ and analytic spread
$\ell$. Let $X_{jl}$, $1 \leq j \leq n$, $1 \leq l \leq \ell$, be
variables, write $S=R(\{ X_{jl}\}) = R[\{ X_{jl} \}]_{{\mathfrak
m}[\{ X_{jl}\}]}$, and consider the $S$-ideal ${\mathcal A}$
generated by $\displaystyle\sum_{j=1}^n X_{jl} f_j$, $1 \leq l
\leq \ell$. This ideal, which we dub a {\it universal
$\ell$-generated ideal in $IS$}, is a minimal reduction of $IS$.
In \cite{RS}, D. Rees and J. Sally prove that if $I$ is
${\mathfrak m}$-primary, then ${\mathcal A} \cap R \subset {\rm
core}(I)$. One of the main results of the present paper says that
this containment is actually an equality $($Theorem~\ref{4.7}.{\it
b}$)$, which remains valid if $R$ is merely Buchsbaum
$($Remark~\ref{4.10}$)$. Most notably however, we are able to
treat ideals that are not necessarily ${\mathfrak m}$-primary,
such as ideals with $\ell=g$ $($called {\it equimultiple}
ideals$)$, generically complete intersection Cohen--Macaulay
ideals in a Gorenstein ring with $\ell = g+1$, or two-dimensional
Cohen--Macaulay ideals in a Gorenstein ring which are complete
intersections locally on the punctured spectrum. All these ideals
fall into the class of {\it universally weakly
$(\ell-1)$-residually $S_2$ ideals satisfying $G_{\ell}$}, a class
that provides the framework for this article $($see Sect. 2 for
the definition$)$. This assumption is fairly general, it
essentially requires the vanishing of $\ell-g$ local cohomology
modules. For such ideals we are able to prove that
\begin{description}[(3)]
\item[$($1$)$]
${\rm core}(I)$ is the intersection of finitely many general minimal
reductions of $I$ $($Theorem~\ref{4.5}$)$;

\item[$($2$)$]
${\rm core}(I) = {\mathcal A} \cap R$ $($Theorem~\ref{4.7}.{\it b}$)$;

\item[$($3$)$]
${\rm core}(IR') = {\rm core}(I)R'$ for every flat $($not necessarily
local$)$ homomorphism $R \longrightarrow R'$ of local Cohen--Macaulay rings
so that $IR'$ is universally weakly $(\ell-1)$-residually $S_2$
$($Theorem~\ref{4.8}$)$.
\end{description}

The main technical result is the fact that ${\rm core}(I)$ can be obtained
by intersecting {\it general} minimal reductions of $I$. It immediately
implies the flat ascent asserted in $($3$)$, provided the map is local
$($Lemma~\ref{4.6}$)$. This yields $($2$)$, which in turn leads to the
general case of $($3$)$.
From the equality in $($2$)$ we also deduce an expression for ${\rm core}(I)$
as a colon ideal in a polynomial ring over $R$ that allows -- at least in
principle -- for an explicit computation of the core $($Proposition~\ref{5.4}
and Remark~\ref{5.5}$)$.

The assertions $($1$)$ and $($2$)$ above are no longer true for arbitrary
ideals in Cohen--Macaulay rings $($Example~\ref{4.11}$)$. On the other hand,
we are able to prove under fairly weak assumptions
that ${\rm core}(I)$ is still
an intersection of {\it finitely} many minimal reductions of $I$, which is
far from being obvious for non ${\mathfrak m}$-primary ideals
$($Theorem~\ref{3.1}$)$. In fact we do not know of any examples where this
finiteness assertion or the flat ascent as in $($3$)$ fail to hold. Thus
we are led to ask the following questions, where $I$ is an arbitrary
$R$-ideal and ${\mathcal M}(I)$ the set of its minimal reductions:
\begin{description}[(iii)]
\item[$(i)$]
Is $\displaystyle\bigcup_{J \in {\mathcal M}(I)} {\rm Ass}_R(I/J)$ finite?

\item[$(ii)$]
Is ${\rm core}(I)$ an intersection of finitely many minimal reductions of $I$?

\item[$(iii)$]
Is ${\rm core}(I) \subset {\mathcal A} \cap R$?

\item[$(iv)$]
Is ${\rm core}(IR') \supset ({\rm core}(I))R'$ for every flat local
homomorphism $R \longrightarrow R'$ of local Cohen--Macaulay rings?
\end{description}
Notice that an affirmative answer to $(iv)$ would imply that
$(iii)$ holds, and that $(ii)$ and $(iv)$, if valid, would yield
the equality ${\rm core}(IR') = ({\rm core}(I))R'$ in the setting
of $(iv)$.

We end by remarking that effective `closed formulas' for the
computation of ${\rm core}(I)$ will appear in another article of ours
\cite{CPU}, extending earlier work by C. Huneke and I. Swanson \cite{HS}.
However, the assumptions on $I$ will be more restrictive than the ones used
here and the techniques will be different.

\section{Definitions and preliminaries}

We begin by reviewing some facts from \cite{CEU} about residually
$S_2$ ideals. Let $R$ be a Noetherian ring, $I$ an $R$-ideal of
height $g$, and $s$ an integer. Recall that $I$ satisfies the
condition $G_s$ if for each prime ideal ${\mathfrak p}$ containing
$I$ with ${\rm dim}\, R_{\mathfrak p} \leq s-1$, the minimal
number of generators $\mu(I_{\mathfrak p})$ is at most ${\rm
dim}\, R_{\mathfrak p}$. A proper $R$-ideal $K$ is called an {\it
$s$-residual intersection} of $I$ if there exists an $s$-generated
ideal $J \subset I$ so that $K= J \colon I$ and ${\rm ht}\, K \geq
s \geq g$. If in addition ${\rm ht}\, I + K \geq s+1$ we say that
$K$ is a {\it geometric $s$-residual intersection} of $I$. The
ideal $I$ is called {\it $s$-residually $S_2$ $($weakly
$s$-residually $S_2)$} if $R/K$ satisfies Serre's condition $S_2$
for every $i$-residual intersection $($geometric $i$-residual
intersection, respectively$)$ $K$ of $I$ and every $i \leq s$.
Finally, whenever $R$ is local, we say $I$ is {\it universally
$s$-residually $S_2$ $($universally weakly $s$-residually $S_2)$}
if $IS$ is residually $S_2$ $($weakly $s$-residually $S_2$,
respectively$)$ for every ring $S=R(X_1, \ldots, X_n)$ with $X_1,
\ldots, X_n$ variables over $R$.

If $(R, {\mathfrak m})$ is a local Cohen--Macaulay ring of dimension $d$
and $I$ an $R$-ideal satisfying $G_s$, then $I$ is universally
$s$-residually $S_2$ in the following cases:
\begin{description}[(2)]
\item[$($1$)$]
$I$ has sliding depth, which means that
the $i^{\rm th}$ Koszul homology modules $H_i$ of a generating set
$f_1, \ldots, f_n$ of $I$ satisfy ${\rm depth}\, H_i \geq d -
n + i$ \, for every $i$ $($see \cite[3.3]{HVV}$)$.

\item[$($2$)$]
$R$ is Gorenstein, and the local cohomology modules $H_{\mathfrak
m}^{d-g-j}(R/I^j)$ vanish for $1 \leq j \leq s-g+1$ or,
equivalently, ${\rm Ext}_R^{g+j}(R/I^j, R) = 0$ for $1 \leq j \leq
s-g+1$ $($see \cite[4.1 and 4.3]{CEU}$)$. The latter condition
holds whenever ${\rm depth}\, R/I^j \geq {\rm dim}\, R/I - j + 1$
for $1 \leq j \leq s-g+1$.
\end{description}

The depth inequalities in $($1$)$ and $($2$)$ are satisfied by
strongly Cohen--Macaulay ideals, i.e., ideals whose Koszul
homology modules are Cohen--Macaulay. The latter condition always
holds if $I$ is a Cohen--Macaulay almost complete intersection or
a Cohen--Macaulay deviation two ideal of a Gorenstein ring
\cite[p.{\,}259]{AH}. It is also satisfied for any ideal in the
linkage class of a complete intersection \cite[1.11]{H1}: standard
examples include perfect ideals of grade $2$ and perfect
Gorenstein ideals of grade $3$.

Finally, $I$ is universally $s$-residually $S_2$ if $s<g$ or if
one of the following conditions holds:
\begin{description}[(3)]
\item[$($3$)$]
$R$ is Gorenstein, $R/I$ is Cohen--Macaulay, and $s=g$ $($see $($2$))$.

\item[$($4$)$]
$R$ is Gorenstein, $R/I$ is Cohen--Macaulay, and ${\rm dim}\, R/I \leq 2$
$($see $($3$)$ and Lemma~\ref{2.2}.{\it b} below$)$.

\end{description}

\begin{lemma}\label{2.2}
Let $R$ be a Cohen--Macaulay ring, $s$
an integer, and $I$ an $R$-ideal satisfying $G_s$.
\begin{description}[(a)]
\item[$(${\it a}$)$]
If $I$ is weakly $(s-1)$-residually $S_2$, then for every ${\mathfrak p}
\in V(I)$, $I_{\mathfrak p}$ is weakly $(s-1)$-residually $S_2$.

\item[$(${\it b}$)$]
If $I$ is weakly $(s-1)$-residually $S_2$, then every $s$-residual
intersection of $I$ is unmixed of height $s$.

\item[$(${\it c}$)$]
If $I$ is weakly $(s-1)$-residually $S_2$, then for every $s$-residual
intersection $J \colon I$ of $I$ with $J \subset I$ and $\mu(J) \leq s$,
every associated prime of $J$ has height at most $s$.

\item[$(${\it d}$)$]
If $I$ is weakly $(s-1)$-residually $S_2$, then for every geometric
$s$-residual intersection $K = J \colon I$ of $I$ with $J \subset I$ and
$\mu(J) \leq s$, $I \cap K =J$.

\item[$(${\it e}$)$]
If $I$ is weakly $(s-2)$-residually $S_2$, then for every flat
homomorphism of Noetherian rings $R \longrightarrow R'$ and every
$s$-generated reduction $J$ of $IR'$, ${\rm Supp}_{R'}(IR'/J) =
V({\rm Fitt}_s(I)R')$.

\item[$(${\it f}$)$]
If $I$ is weakly $(s-1)$-residually $S_2$, then for every $s$-generated
reduction $J$ of $I$, ${\rm Ass}_R(I/J)$ is empty or consists only of
primes of height $s$.

\item[$(${\it g}$)$]
If $I$ is weakly $(s-2)$-residually $S_2$, then $I=J$ for every
$(s-1)$-generated reduction $J$ of $I$.
\end{description}
\end{lemma}
\begin{proof}
Part $(${\it a}$)$ follows as in \cite[the proof of 1.10.{\it a}]{AN}, whereas
parts $(${\it b}$)$, $(${\it c}$)$ and $(${\it d}$)$ are identical
to \cite[3.4.{\it a},\,{\it b},\,{\it c}]{CEU}.

To prove part $(${\it e}$)$ it suffices to show that if ${\mathfrak q}
\in V(IR')$ and $\mu(IR_{\mathfrak q}') \leq s$, then $IR_{\mathfrak q}'
= J_{\mathfrak q}'$. So let ${\mathfrak p} = {\mathfrak q} \cap R$. Now
$\mu(I_{\mathfrak p}) \leq s$ and by part $(${\it a}$)$,
$I_{\mathfrak p}$ is still weakly $(s-2)$-residually $S_2$. Thus
according to \cite[3.6.{\it b}]{CEU}, $I_{\mathfrak p}$ can be generated
by a $d$-sequence. Therefore $\mu(I_{\mathfrak p})=\ell(I_{\mathfrak p})$
by \cite[3.1]{Huneke} or \cite[3.15]{Valla},
and hence $\mu(IR_{\mathfrak q}') =\ell(IR_{\mathfrak
q}')$. Since $J_{\mathfrak q}$ is a reduction of $IR_{\mathfrak q}'$, we
conclude that $IR_{\mathfrak q}' =J_{\mathfrak q}$.

As to part $(${\it f}$)$ we may assume that $I \not=J$. By $(${\it e}$)$,
${\rm codim}\, {\rm Supp}_R(I/J) \geq s$.
In particular $J \colon I$ is an $s$-residual intersection of $I$, which
according to $(${\it c}$)$ implies that every prime in ${\rm
Ass}_R(R/J)$ has height at most $s$. Now our assertion follows since
${\rm Ass}_R(I/J) \subset {\rm Supp}_R(I/J) \cap {\rm Ass}_R(R/J)$.

Finally we show part $(${\it g}$)$. Applying $(${\it e}$)$ with
$s$ replaced by $s-1$, we deduce that ${\rm codim}\, {\rm
Supp}_R(I/J) > s-1$, which implies $I=J$ by $(${\it f}$)$. \QED
\end{proof}

We are now going to introduce the notions of generic, universal, and
general subideals, that will play a crucial role in the sequel.
To this end, let $R$ be a Noetherian ring, $I$ an $R$-ideal, $f_1, \ldots,
f_n$ a generating sequence of $I$, and $t$, $s$ integers. Let
$\underline{\underline{X}} = X_{ijl}$,
$1 \leq i \leq t$, $1 \leq j \leq n$, $1 \leq l \leq s$, be variables,
$T=R[\underline{\underline{X}}] =
R[\{ X_{ijl}\}]$, and ${\mathcal B}_i$, $1 \leq i \leq t$, the $T$-ideal
generated by $\displaystyle\sum_{j=1}^n X_{ijl} f_j$, $1 \leq l \leq s$.
We call ${\mathcal B}_1, \ldots, {\mathcal B}_t$ {\it generic $s$-generated
ideals in $IT$}.
Notice that up to adjoining variables and applying an $R$-automorphism, this
definition does not depend on the choice of $n$ and $f_1, \ldots, f_n$.

Now assume $(R, {\mathfrak m})$ is local with residue field $k$. We set
$S = T_{{\mathfrak m}T}$ and ${\mathcal A}_i = {\mathcal B}_iS$, and call
${\mathcal A}_1, \ldots, {\mathcal A}_t$ {\it universal $s$-generated
ideals in $IS$}.
Furthermore we say that ${\mathfrak a}_1, \ldots, {\mathfrak a}_t$ are
{\it general $s$-generated ideals in $I$} if ${\mathfrak a}_i \subset I$
are ideals with $\mu({\mathfrak a}_i)=s$, $ {\mathfrak a}_i \otimes_R
k \hookrightarrow I \otimes_R k$, and the point $({\mathfrak
a}_1\otimes_R k, \ldots, {\mathfrak a}_t \otimes_R k)$ lies in
some dense open subset of the product of Grassmannians
$\displaystyle\stackrel{t}{\times} {\mathbb G}(s, I \otimes_R k)$.

Now let $\underline{\underline{\lambda}} = \lambda_{ijl}$, $1 \leq i
\leq t$, $1 \leq j \leq n$,
$1 \leq l \leq s$ be elements of $R$, and consider the maximal ideal
$M = ({\mathfrak m}, \underline{\underline{X}}-\underline{\underline{\lambda}}) =
({\mathfrak m}, \{ X_{ijl}-\lambda_{ijl} \})$ of $T$.
We will identify the set $\{ M = ({\mathfrak m},
\underline{\underline{X}}-\underline{\underline{\lambda}}) \,|\,
\underline{\underline{\lambda}} \in R^{t n s} \}$ with the set of
$k$-rational points of the affine space ${\mathbb A}_k^{t n s}$. Write $\pi =
\pi_{\underline{\underline{\lambda}}} : T \longrightarrow
R$ for the homomorphism of $R$-algebras with $\pi(X_{ijl}) = \lambda_{ijl}$.
The kernel of $\pi$ is generated by the $T$-regular sequence
$\underline{\underline{X}}-\underline{\underline{\lambda}}$.
Now $\pi({\mathcal B}_1), \ldots, \pi({\mathcal B}_t)$ is a sequence of
$s$-generated ideals in $I$, whose images in $I \otimes_R k$ only depend
on $M$. Conversely, every sequence of $s$-generated ideals in $I$ is
obtained in this way. As $\underline{\underline{X}}-
\underline{\underline{\lambda}}$ form a regular sequence modulo every
power of $IT$, Nakayama's Lemma shows that $\pi({\mathcal B}_i)$ is a
reduction of $I$ if and only if $({\mathcal B}_i)_M$ is a reduction
of $IT_M$. Since the latter condition is equivalent to
$M \not\in \displaystyle \bigcap_{r \geq 0} {\rm Supp}_T (I^{r+1}T/{\mathcal
B}_iI^r)$, the set of all $M$ for which $({\mathcal B}_i)_M$ is a reduction
of $IT_M$ is open in ${\mathbb A}_k^{t n s}$.

Finally, we write ${\mathcal M}(I)$ for the set of all
minimal reductions of $I$,  and we define
\[
{\mathcal P}(I) = \bigcup_{J \in {\mathcal M}(I)} {\rm Ass}_R(I/J)
\]
and
\[
\gamma(I) = {\rm inf}\, \{ t \, | \, {\rm core}(I) = \bigcap_{i=1}^t J_i \
{\rm with} \ J_i \in {\mathcal M}(I) \}.
\]

\section{Finiteness}

If $I$ is ${\mathfrak m}$-primary,
then ${\rm core}(I)$ is ${\mathfrak m}$-primary and $\gamma(I) \leq
{\rm type}\, (R/{\rm core}(I))$.
In particular $\gamma(I)$ is finite. The next result establishes this
finiteness in a much broader setting.

\begin{theorem}\label{3.1}
Let $R$ be a Noetherian local ring with infinite residue field and $I$
an $R$-ideal. Assume that ${\mathcal P}(I)$ is finite, every element
of ${\mathcal P}(I)$ is minimal in this set, and $\mu(I_{\mathfrak q})
= \ell(I_{\mathfrak q})$ for every ${\mathfrak q} \in {\rm Spec}(R)
\setminus \overline{{\mathcal P}(I)}$. Then $\gamma(I) < \infty$.
\end{theorem}
\begin{proof}
Let $K$ be the intersection of the prime ideals in
${\mathcal P} = {\mathcal P}(I)$. We first prove that
$K \subset \sqrt{{\rm core}(I) \colon I}$, or equivalently that there exists
a fixed integer $r$ so that $K^r I \subset J$ for every
$J \in {\mathcal M} = {\mathcal M}(I)$.

Let ${\mathfrak m}$ be the maximal ideal of $R$, $k = R/{\mathfrak m}$,
$\ell = \ell(I)$, ${\mathcal B}$ a generic $\ell$-generated ideal in $IT$,
and $U \subset {\mathbb A}_k^{n \ell}$ the open set consisting of all $M =
({\mathfrak m}, \underline{\underline{X}}-\underline{\underline{\lambda}})
\subset T$ such that ${\mathcal B}_M$ is a reduction of $IT_M$. Fix $M \in
U$ and let $Q$ be any prime ideal of $T$ with $K \not\subset
Q \subset M$. Writing ${\mathfrak q} = Q \cap R$ we have ${\mathfrak q}
\in {\rm Spec}(R) \setminus \overline{{\mathcal P}}$. Thus
$\mu(I_{\mathfrak q}) = \ell(I_{\mathfrak q})$ and hence $\mu(IT_Q) =
\ell(IT_Q)$. Since ${\mathcal B}_Q$ is a reduction of $IT_Q$ we conclude that
$IT_Q = {\mathcal B}_Q$. Therefore $K^n IT_M \subset {\mathcal B}_M$ for some
integer $n$. Now $K^nIT_N \subset {\mathcal B}_N$ for every $N$ in some
open neighborhood of $M$ in $U$. As $U$ is quasicompact there exists
an integer $r$ so that $K^r IT_N \subset {\mathcal B}_N$ for every
$N \in U$. Specializing we conclude that $K^rI \subset J$ for every
$J \in {\mathcal M}$.

Now for every ${\mathfrak p} \in {\mathcal P}$, since ${\mathfrak
p}$ is minimal in the finite set ${\mathcal P}$, ${\mathfrak p}^r
I_{\mathfrak p} = K^r I_{\mathfrak p} \subset ({\rm
core}(I))_{\mathfrak p}$ and hence ${\rm length}_{R_{\mathfrak
p}}((I/{\rm core}(I))_{\mathfrak p}) < \infty$. Again as
${\mathcal P}$ is finite there exist finitely many minimal
reductions $J_1, \ldots, J_t$ so that for every ${\mathfrak p} \in
{\mathcal P}$, $\displaystyle \bigcap_{i=1}^t (J_i/{\rm
core}(I))_{\mathfrak p} = \displaystyle \bigcap_{J \in {\mathcal
M}} (J/{\rm core}(I))_{\mathfrak p}$. Thus $\displaystyle
(\bigcap_{i=1}^t J_i)_{\mathfrak p} \subset J_{\mathfrak p}$ for
every $J \in {\mathcal M}$ and every ${\mathfrak p} \in {\mathcal
P}$. Hence by the definition of ${\mathcal P}$, $\displaystyle
\bigcap_{i=1}^t J_i \subset J$, which gives ${\rm core}(I) =
\displaystyle \bigcap_{i=1}^t J_i$. \QED \end{proof}

\begin{remark}\label{3.2}
The assumptions on ${\mathcal P}(I)$ in Theorem~\ref{3.1}
are automatically satisfied in any of the following cases, where
$g = {\rm ht}(I)$ and $\ell = \ell(I)$:

\begin{description}[$\bullet$]
\item[$\bullet$]
Locally on the punctured spectrum of $R$, $I$ is generated by
analytically independent elements.

\item[$\bullet$]
$R$ satisfies $S_{g+1}$ and $I$ is equimultiple.

\item[$\bullet$]
$R$ is Cohen--Macaulay and $I$ is $G_{\ell}$ and weakly
$(\ell-1)$-residually $S_2$ $($see Lemma~\ref{2.2}.{\it a},\,{\it
e},\,{\it f}$)$.
\end{description}

In either case ${\mathcal P}(I) = {\rm Min}({\rm Fitt}_{\ell}(I)) =
{\rm Ass}_R(I/J)$, where $J$ is any minimal reduction of $I$.
\end{remark}

\section{Genericity and the shape of the core}

The crucial result of this section is Theorem~\ref{4.5}, which describes
the core as a finite intersection of {\it general} minimal reductions.
To prove it we compare the multiplicities of modules defined by intersecting
reduction ideals, universal ideals, and general ideals, respectively.
This is done in Lemma~\ref{4.4}, which in turn is based on Lemmas~\ref{4.1}
and~\ref{4.2}. There we deal with modules defined by a single universal ideal
$($Lemma~\ref{4.1}$)$ and a single minimal reduction $($Lemma~\ref{4.2}$)$,
respectively. For the latter we prove that under suitable assumptions, the
multiplicity of $I/J$ is independent of the choice of a minimal reduction $J$
of $I$, a fact reminiscent of the theme of \cite{CEU}.

If $R$ is a Noetherian local ring and $E$ a finitely generated $R$-module, we
denote by $e(E)$ the multiplicity of $E$, by $e_I(E)$ the multiplicity with
respect to an ideal of definition $I$ of $E$, and by $e(\underline{y}; E)$ the
multiplicity with respect to a system of parameters $\underline{y}$ of $E$.

\begin{lemma}\label{4.1}
Let $R$ be a Noetherian local ring, $I$ an $R$-ideal, and $s$ an
integer. Let $J \subset I$ be an $s$-generated ideal and
${\mathcal A}$ a universal $s$-generated ideal in $IS$. One has
${\rm dim}\, I/J \geq {\rm dim}\, IS/{\mathcal A}$, and if
equality holds then $e(I/J) \geq e(IS/{\mathcal A})$.
\end{lemma}
\begin{proof}
We may assume that the residue field of $R$ is infinite. Let
${\mathfrak m}$ be the maximal ideal of $R$, and ${\mathcal B}$ a
generic $s$-generated ideal in $IT$ so that ${\mathcal A} =
{\mathcal B}S ={\mathcal B}_{{\mathfrak m}T}$. Using the notation
introduced in Sect. 2 we have $J =
\pi_{\underline{\underline{\lambda}}}({\mathcal B})$ for some
$\underline{\underline{\lambda}} \in R^{ns}$. Write $M =
({\mathfrak m}, \underline{\underline{X}}-
\underline{\underline{\lambda}}) \subset T$. First notice that
$I/J \simeq (IT/{\mathcal B}+(\underline{\underline{X}}-
\underline{\underline{\lambda}})IT)_M$ as
$\underline{\underline{X}}-\underline{\underline{\lambda}}$ form a
regular sequence on $(T/IT)_M$. In particular ${\rm dim}\, I/J
\geq {\rm dim}\, (IT/{\mathcal B})_M -ns$. On the other hand ${\rm
dim}\, (IT/{\mathcal B})_M - ns = {\rm dim}\, (IT/{\mathcal B})_M
- {\rm dim}\, (T/{\mathfrak m}T)_M \geq {\rm dim}\, (IT/{\mathcal
B})_{{\mathfrak m}T}$. Now the inequality ${\rm dim}\, I/J \geq
{\rm dim}\, IS/{\mathcal A}$ follows.

Moreover if equality holds then ${\rm dim}\, I/J = {\rm dim}\,
(IT/{\mathcal B})_M - ns = {\rm dim}\, (IT/{\mathcal
B})_{{\mathfrak m}T}$. Hence there exists a sequence
$\underline{y}$ in ${\mathfrak m}$ so that $\underline{y}$ is a
system of parameters of $I/J$ and of $(IT/{\mathcal
B})_{{\mathfrak m}T}$, $\underline{y}$ generates a minimal
reduction of ${\mathfrak m}/J \colon I$, and $\underline{y},
\underline{\underline{X}}-\underline{\underline{\lambda}}$ form a
system of parameters of $(IT/{\mathcal B})_M$. Therefore $e(I/J) =
e(\underline{y}; I/J) \geq e(\underline{y},
\underline{\underline{X}}- \underline{\underline{\lambda}};
(IT/{\mathcal B})_M) \geq e(\underline{y}; (IT/{\mathcal
B})_{{\mathfrak m}T})$ where the last inequality holds by the
associativity formula for multiplicities as ${\rm dim}\,
(T/{\mathfrak m}T)_M = ns$ $($see \cite[24.7]{Nagata}$)$. Finally
$e(\underline{y}; (IT/{\mathcal B})_{{\mathfrak m}T}) \geq
e(IS/{\mathcal A})$. \QED \end{proof}

\begin{lemma}\label{4.2}
Let $R$ be a local Cohen--Macaulay ring with infinite residue field and
$I$ an $R$-ideal of analytic spread $\ell$. Assume that $I$ satisfies
$G_{\ell}$, but not $G_{\ell+1}$, and $I$ is weakly $(\ell-2)$-residually
$S_2$. Then $e(I/J)$ and $e(R/J \colon I)$ are independent of the minimal
reduction $J$ of $I$.
\end{lemma}
\begin{proof}
Let $J$ be any minimal reduction of $I$. By Lemma~\ref{2.2}.{\it
e}, ${\rm Supp}(I/J) = V({\rm Fitt}_{\ell}(I))$, and by our
assumption the latter set has codimension $\ell$. According to the
associativity formula for multiplicities we may localize at any
minimal prime in ${\rm Supp}(I/J) = V({\rm Fitt}_{\ell}(I))$ of
codimension $\ell$ to assume that ${\rm dim}\, R = \ell$. By
Lemma~\ref{2.2}.{\it a}, $I$ is still weakly $(\ell-2)$-residually
$S_2$. Notice that now ${\rm dim}\, I/J = {\rm dim}\, R/J \colon I
= 0$.

We may suppose $\ell > 0$. Let ${\mathfrak a}$ and ${\mathfrak b}$
be minimal reductions of $I$. By a general position argument
$($see, e.g., \cite[the proof of 1.4]{AN}$)$ there are generating
sequences $a_1, \ldots, a_{\ell}$ of ${\mathfrak a}$ and $b_1,
\ldots, b_{\ell}$ of ${\mathfrak b}$ so that for every $1 \leq i
\leq \ell$ the ideal $(a_1, \ldots, a_i, b_{i+1},\ldots,
b_{\ell})$ is a minimal reduction of $I$ and
\[
K =(a_1, \ldots, a_{i-1}, b_{i+1}, \ldots, b_{\ell}) \colon I
\]
is a geometric $(\ell-1)$-residual intersection of $I$.
Write ${\mathfrak c} = (a_1, \ldots, a_{i-1}, b_{i+1},$ $\ldots, b_{\ell})$.
It suffices to prove that
\[
e(I/({\mathfrak c}, a_i)) = e(I/({\mathfrak c}, b_i))
\qquad {\rm and} \qquad
e(R/({\mathfrak c}, a_i) \colon I) = e(R/({\mathfrak c}, b_i) \colon I).
\]

To this end let `$\, {}^{\tratto}$' denote images in
$\overline{R}=R/K$. Since $I \cap K = {\mathfrak c}$
by Lemma~\ref{2.2}.{\it d}, it follows that
$\overline{I}/(\overline{a}_i) = I/(I \cap K, a_i) =
I/({\mathfrak c}, a_i)$ and hence
$(\overline{a}_i) \colon \overline{I} =
\overline{({\mathfrak c}, a_i) \colon I}$.
The same holds for $b_i$ in place of $a_i$. Thus it suffices to show
$e(\overline{I}/(\overline{a}_i)) = e(\overline{I}/(\overline{b}_i))$
and $e(\overline{R}/(\overline{a}_i) \colon \overline{I}) =
e(\overline{R}/(\overline{b}_i) \colon \overline{I})$. Obviously,
$(\overline{a}_i)$ and $(\overline{b}_i)$ are reductions of
$\overline{I}$, and by Lemma~\ref{2.2}.{\it b}, $\overline{a}_i$
and $\overline{b}_i$ are non zerodivisors on $\overline{R}$ and
${\rm dim}\, \overline{R} = 1$.

After changing notation, we are reduced to proving that if
$(R, {\mathfrak m})$
is a one-dimensional local Cohen--Macaulay ring with infinite residue
field and $I$ is an ${\mathfrak m}$-primary ideal, then
${\rm length}(I/J)$ and ${\rm length}(R/J \colon I)$ do not depend on the
minimal reduction $J$ of $I$. Notice that ${\rm length}(R/J) = e_J(R) =
e_I(R)$. This gives ${\rm length}(I/J) = e_I(R)-{\rm length}(R/I)$.
If $a$ is an $R$-regular element generating $J$, we have
$I^{-1}/R \simeq aI^{-1}/(a)$ and $aI^{-1} = J \colon I$.
This yields an  exact sequence
\[
0 \rightarrow I^{-1}/R \longrightarrow R/J
\longrightarrow R/J \colon I \rightarrow 0,
\]
which shows that ${\rm length}(R/J \colon I) = e_I(R) - {\rm
length}(I^{-1}/R)$. \QED \end{proof}

\begin{lemma}\label{4.4}
Let $R$ be a local Cohen--Macaulay ring with infinite residue field,
$I$ an $R$-ideal of analytic spread $\ell$, and $J_1, \ldots, J_t$
minimal reductions of $I$. Assume that $I$ is $G_{\ell}$, but not
$G_{\ell+1}$, and $I$ is weakly $(\ell-2)$-residually $S_2$.
Then
\[
e(I/J_1 \cap \ldots \cap J_t) \leq e(IS/{\mathcal
A}_1 \cap \ldots \cap {\mathcal A}_t)= e(I/{\mathfrak a}_1
\cap \ldots \cap {\mathfrak a}_t)
\]
and ${\mathfrak a}_1, \ldots, {\mathfrak a}_t$ are minimal
reductions of $I$, for ${\mathcal A}_1, \ldots, {\mathcal A}_t$
universal $\ell$- \!\!\! generated ideals in $IS$ and ${\mathfrak a}_1,
\ldots, {\mathfrak a}_t$ general $\ell$-generated ideals in $I$.
\end{lemma}
\begin{proof}
We may assume that $n=\mu(I) > \ell > 0$. Write $d = {\rm dim}\,
(R/{\rm Fitt}_{\ell}(I))$. By Lemma~\ref{2.2}.{\it e},
\begin{equation}\label{eq1}
{\rm Supp}_{R'}(IR'/J) = V({\rm Fitt}_{\ell}(I)R')
\end{equation}
for any flat homomorphism of rings $R \longrightarrow R'$ and any
$\ell$-generated reduction $J$ of $IR'$.
Let ${\mathcal B}_1, \ldots, {\mathcal B}_t$ be generic $\ell$-generated
ideals in $IT$, defined using a minimal generating sequence of $I$. We may
suppose that ${\mathcal A}_i = {\mathcal B}_i S$.

We first prove the equality $e(IS/{\mathcal A}_1 \cap \ldots \cap
{\mathcal A}_t)
= e(I/ {\mathfrak a}_1 \cap \ldots \cap {\mathfrak a}_t)$. We may pass from
$R$ to $\widehat{R}$ and assume that $(R, {\mathfrak m}, k)$ is complete.
Although the residually $S_2$ assumption may not be preserved by completion,
Lemma~\ref{2.2}.{\it e} shows that $($\ref{eq1}$)$ still holds.
Let ${\mathfrak p}_1, \ldots, {\mathfrak p}_s$ be the minimal primes of
${\rm Fitt}_{\ell}(I)$ having maximal dimension, namely $d$.
Let $C$ be a coefficient ring of $R$ and
write $C_j = C/{\mathfrak p}_j \cap C$, $k_j = {\rm Quot}(C_j)$ for $1
\leq j \leq s$.

For $1 \leq i \leq t$ consider the $T$-modules $E_i = IT/{\mathcal
B}_1 \cap \ldots \cap {\mathcal B}_i + {\mathcal B}_{i+1}$, where
${\mathcal B}_{t+1} = 0$. One has ${\rm dim}(E_i)_{{\mathfrak
p}_jT} \leq 0$ by $($\ref{eq1}$)$. Thus the set of all maximal
ideals $({\mathfrak p}_j, \underline{\underline{X}}-
\underline{\underline{\lambda}})$ of $T_{{\mathfrak p}_j}$ with
$\underline{\underline{\lambda}} \in C_{{\mathfrak p}_j \cap C}^{t
n \ell}$, so that the modules $(E_i)_{({\mathfrak p}_j,
\underline{\underline{X}}- \underline{\underline{\lambda}})}$ are
zero or Cohen--Macaulay of dimension $tn\ell$ for $1 \leq i \leq
t$, form a dense open subset of ${\mathbb A}_{k_j}^{t n \ell}$.
Let $H_j$ be the largest ideal of $k_j[\underline{\underline{X}}]$
defining the complement of this subset, and let $\overline{H}_j$
denote the image of $H_j \cap C_j[\underline{\underline{X}}]$ in
$k[\underline{\underline{X}}]$. Since $C_j$ is either a field or a
discrete valuation ring, it follows that $\overline{H}_j \not= 0$
and hence $U_j = D(\overline{H}_j)$ is a dense open subset of
${\mathbb A}_k^{t n \ell}$. Let $M$ be a maximal ideal in $T$ of
the form $M = ({\mathfrak m},
\underline{\underline{X}}-\underline{\underline{\lambda}})$ with
$\underline{\underline{\lambda}} \in R^{t n \ell}$. As $C$ is a
coefficient ring of $R$ we may assume that
$\underline{\underline{\lambda}} \in C^{t n \ell}$. If the image
of $\underline{\underline{\lambda}}$ in ${\mathbb A}_k^{t n \ell}$
lies in $U_j = D(\overline{H}_j)$, then the image of
$\underline{\underline{\lambda}}$ in ${\mathbb A}_{k_j}^{t n
\ell}$ belongs to $D(H_j)$ since $C_j$ maps onto $k$. Thus the
modules $(E_i)_{ ({\mathfrak p}_j,
\underline{\underline{X}}-\underline{\underline{\lambda}})}$ are
zero or Cohen--Macaulay of dimension $tn\ell$ for $1 \leq i \leq
t$, whenever $M = ({\mathfrak m},
\underline{\underline{X}}-\underline{\underline{\lambda}})$ lies
in $U_j$. Finally, let $U_0$ be the dense open subset of ${\mathbb
A}_k^{t n \ell}$ consisting of all $M=({\mathfrak m},
\underline{\underline{X}}-\underline{\underline{\lambda}})$ so
that $\pi_{\underline{\underline{\lambda}}}({\mathcal B}_i)$ are
reductions of $I$ for $1 \leq i \leq t$. Define $U$ to be the
dense open subset $U_0 \cap \ldots \cap U_s$ of ${\mathbb A}_k^{t
n \ell}$. The natural action of $\displaystyle\stackrel{t}{\times}
G\!L_{\ell}(k)$ on ${\mathbb A}_k^{t n \ell} = \displaystyle
\stackrel{t}{\times} {\rm Hom}_k(k^{\ell},k^n)$ induces an action
on $U$, and so the image $V$ of $U$ in the product of
Grassmannians $\displaystyle\stackrel{t}{\times} {\mathbb G}(\ell,
k^n)$ is open and dense. It remains to show that $e(IS/{\mathcal
A}_1 \cap \ldots \cap {\mathcal A}_t) = e(I/{\mathfrak a}_1 \cap
\ldots \cap {\mathfrak a}_t)$ whenever $({\mathfrak a}_1 \otimes_R
k, \ldots, {\mathfrak a}_t \otimes_R k) \in V$.

So let $({\mathfrak a}_1 \otimes_R k, \ldots, {\mathfrak a}_t
\otimes_R k) \in V$. Write ${\mathfrak a}_i =
\pi_{\underline{\underline{\lambda}}} ({\mathcal B}_i)$ where
$\underline{\underline{\lambda}} \in C^{t n \ell}$, $M=({\mathfrak
m}, \underline{\underline{X}}- \underline{\underline{\lambda}})T$,
and $Q_j=({\mathfrak p}_j,
\underline{\underline{X}}-\underline{\underline{\lambda}})T$. By
the above ${\mathfrak a}_i$ are reductions of $I$, $({\mathcal
B}_i)_M$ are reductions of $IT_M$, and $(E_i)_{Q_j}$ are zero or
Cohen--Macaulay of dimension $tn\ell$ for $1 \leq i \leq t$, $1
\leq j \leq s$. Since the modules $(E_i)_{Q_j}$ are annihilated by
a power of ${\mathfrak p}_j$ according to $($\ref{eq1}$)$, it then
follows that these modules vanish or that
$\underline{\underline{X}}-\underline{\underline{\lambda}}$ form a
regular sequence on them. Notice that $(E_t)_{Q_j} \not= 0$ by
$($\ref{eq1}$)$. Now an induction on $i$, $1 \leq i \leq t$ yields
$(I/{\mathfrak a}_1 \cap \ldots \cap {\mathfrak a}_t)_{{\mathfrak
p}_j} \simeq (IT/{\mathcal B}_1 \cap \ldots \cap {\mathcal B}_t +
(\underline{\underline{X}}-\underline{\underline{\lambda}})IT)_{Q_j}$.
Furthermore, the latter module has length equal to
$e(\underline{\underline{X}}-\underline{\underline{\lambda}};
(IT/{\mathcal B}_1 \cap \ldots \cap {\mathcal B}_t)_{Q_j})$. Since
${\rm Supp}(I/{\mathfrak a}_1 \cap \ldots \cap {\mathfrak a}_t) =
V({\rm Fitt}_{\ell}(I))$ by $($\ref{eq1}$)$, the associativity
formula for multiplicities then yields
\begin{eqnarray*}
e(I/{\mathfrak a}_1 \cap \ldots \cap {\mathfrak a}_t) & = &
\sum_{j=1}^s \lambda_{R_{{\mathfrak p}_j}}((I/{\mathfrak a}_1 \cap
\ldots \cap {\mathfrak a}_t)_{{\mathfrak p}_j}) \cdot e (R/{\mathfrak p}_j) \\
& = & \sum_{j=1}^s
e(\underline{\underline{X}}-\underline{\underline{\lambda}};
(IT/{\mathcal B}_1 \cap \ldots \cap {\mathcal B}_t)_{Q_j}) \cdot
e((T/Q_j)_M).
\end{eqnarray*}

To further evaluate this sum, write $E = E_t = IT/{\mathcal B}_1
\cap \ldots \cap {\mathcal B}_t$. By $($\ref{eq1}$)$, ${\rm
Supp}(E_M) = V({\rm Fitt}_{\ell}(I)T_M)$, which has dimension $d+t
n \ell$. For $r$ an integer with ${\rm Fitt}_{\ell}(I)^r \subset
{\rm ann}(E_M)$, let $\underline{y} = y_1, \ldots, y_d$ be a
sequence of elements generating a minimal reduction of ${\mathfrak
m}/{\rm Fitt}_{\ell}(I)^r$. Now $\underline{y},
\underline{\underline{X}}-\underline{\underline{\lambda}}$ form a
system of parameters of $E_M$. Furthermore in the ring $(T/{\rm
ann}(E))_M$,
$\underline{\underline{X}}-\underline{\underline{\lambda}}$
generate an ideal of height $t n \ell$ and dimension $d$, and
$\underline{y}$ generate an ideal of height $d$ and dimension $t n
\ell$. The minimal primes of maximal dimension of the first ideal
are $Q_1T_M, \ldots, Q_sT_M$, whereas the second ideal has only
one minimal prime, ${\mathfrak m}T_M$. Thus the associativity
formula $($see \cite[24.7]{Nagata}$)$ yields
\begin{eqnarray*}
\sum_{j=1}^s
e(\underline{\underline{X}}-\underline{\underline{\lambda}};
E_{Q_j}) \hspace{-.1cm} & \cdot &  \hspace{-.1cm} e((T/Q_j)_M) =
\sum_{j=1}^s
e(\underline{\underline{X}}-\underline{\underline{\lambda}};
E_{Q_j}) \cdot e(\underline{y}; (T/Q_j)_M) \\
& = & e(\underline{y}, \underline{\underline{X}}-
\underline{\underline{\lambda}}; E_M) = e(\underline{y};
E_{{\mathfrak m}T}) \cdot e(\underline{\underline{X}}-
\underline{\underline{\lambda}}; (T/{\mathfrak m}T)_M)  \\
& = & e(E_{{\mathfrak m}T}) \cdot 1 = e(IS/{\mathcal A}_1 \cap \ldots
\cap {\mathcal A}_t).
\end{eqnarray*}
This completes the proof of the equality $e(IS/{\mathcal A}_1 \cap \ldots
\cap {\mathcal A}_t) = e(I/{\mathfrak a}_1 \cap \ldots
\cap {\mathfrak a}_t)$.

We now show the inequality $e(I/J_1 \cap \ldots \cap J_t) \leq
e(IS/{\mathcal A}_1 \cap \ldots \cap {\mathcal A}_t)$. Writing
${\mathcal J} = J_1S \cap \ldots \cap J_{i-1}S \cap {\mathcal A}_{i+1}
\cap \ldots \cap {\mathcal A}_t$ with $1 \leq i \leq t$, it suffices to
show that $e(IS/{\mathcal J} \cap J_iS) \leq e(IS/{\mathcal J} \cap
{\mathcal A}_i)$.
Consider the two exact sequences
\[
0 \rightarrow IS/{\mathcal J} \cap J_iS \longrightarrow
IS/{\mathcal J} \oplus IS/J_iS \longrightarrow IS/{\mathcal J} + J_iS
\rightarrow 0
\]
\[
0 \rightarrow IS/{\mathcal J} \cap {\mathcal A}_i \longrightarrow
IS/{\mathcal J} \oplus IS/{\mathcal A}_i \longrightarrow IS/{\mathcal J} +
{\mathcal A}_i
\rightarrow 0.
\]
By Lemma~\ref{4.1}, ${\rm dim}\, IS/{\mathcal J}+J_iS \geq {\rm
dim}\, IS/{\mathcal J}+{\mathcal A}_i$, and by $($\ref{eq1}$)$,
all the other $S$-modules occurring in the two exact sequences
have the same dimension, namely $d$. Applying the equality just
proved with $t=1$, gives $e(IS/{\mathcal A}_i) = e(I/{\mathfrak
a}_i)$, whereas Lemma~\ref{4.2} yields $e(I/{\mathfrak a}_i) =
e(I/J_i)$. Thus $e(IS/{\mathcal A}_i) = e(IS/J_iS)$. Now if $d >
{\rm dim}\, IS/{\mathcal J}+J_iS$ or if ${\rm dim}\, IS/{\mathcal
J}+J_iS > {\rm dim}\, IS/{\mathcal J}+{\mathcal A}_i$, then the
asserted inequality follows from the above exact sequences. If on
the other hand $d = {\rm dim}\, IS/{\mathcal J}+J_iS = {\rm dim}\,
IS/{\mathcal J}+{\mathcal A}_i$ then $e(IS/ {\mathcal J} + J_iS)
\geq e(IS/{\mathcal J} + {\mathcal A}_i)$ by Lemma~\ref{4.1}, and
again our claim can be deduced from the two exact sequences. \QED
\end{proof}

\begin{remark}\label{4.3}
The assumption of $I$ not satisfying $G_{\ell+1}$ in
Lemmas~\ref{4.2} and~\ref{4.4} can be omitted if $I$ is weakly
$(\ell-1)$-residually $S_2$; for otherwise $J=I$, $J_i =
{\mathfrak a}_i = I$ and ${\mathcal A}_i = IS$ by
Lemma~\ref{2.2}.{\it e},{\it g}.
\end{remark}

\begin{theorem}\label{4.5}
Let $R$ be a local Cohen--Macaulay ring with infinite residue field
and $I$ an $R$-ideal of analytic spread $\ell$.
Assume that $I$ is $G_{\ell}$ and weakly $(\ell-1)$-residually $S_2$,
and write $t=\gamma(I)$. Then ${\rm core}(I)={\mathfrak a}_1 \cap
\ldots \cap {\mathfrak a}_t$ for ${\mathfrak a}_1, \ldots,
{\mathfrak a}_t$ general {$\ell$-}generated ideals in $I$ which are
reductions of $I$.
\end{theorem}
\begin{proof}
According to Theorem~\ref{3.1} and Remark~\ref{3.2}, $\gamma(I) <
\infty$. Hence by Lemma~\ref{2.2}.{\it f} every associated prime
of the $R$-module $I/{\rm core}(I)$ has height $\ell$, and the
same holds for $I/{\mathfrak a}_1 \cap \ldots \cap {\mathfrak
a}_t$. On the other hand Lemma~\ref{4.4} and Remark~\ref{4.3} show
that $e(I/{\rm core}(I)) \leq e(I/{\mathfrak a}_1 \cap \ldots \cap
{\mathfrak a}_t)$. Now as $I/{\rm core}(I)$ maps onto
$I/{\mathfrak a}_1 \cap \ldots \cap {\mathfrak a}_t$, these two
modules are equal. \QED \end{proof}

Using Theorem~\ref{4.5}, we can now prove that the formation of the core
is compatible with flat local maps:

\begin{lemma}\label{4.6}
Let $R \hookrightarrow R'$ be a flat local extension of local rings
with infinite residue fields. Assume $R'$ is Cohen--Macaulay.
Let $I$ be an $R$-ideal of analytic spread $\ell$ such that $IR'$
is $G_{\ell}$ and weakly $(\ell-1)$-residually $S_2$.
Then ${\rm core}(IR')=({\rm core}(I))R'$ and $\gamma(IR') = \gamma(I)$.
\end{lemma}
\begin{proof}
Notice that $R$ is Cohen--Macaulay and $I$ is $G_{\ell}$ and weakly
$(\ell-1)$-residually $S_2$. By Theorem~\ref{3.1}, ${\rm core}(I) =
\displaystyle \bigcap_{i=1}^s J_i$ for finitely many reductions
$J_1, \ldots, J_s$ of $I$. Thus $({\rm core}(I))R' =
\displaystyle \bigcap_{i=1}^s (J_iR')$. As $J_iR'$ are reductions of $IR'$,
it follows that $({\rm core}(I))R' \supset {\rm core}(IR')$.

To show the other inclusion, let $k \subset K$ be the residue
field extension of $R \hookrightarrow R'$. By Theorem~\ref{4.5},
${\rm core}(IR')$ is the intersection of $t = \gamma(IR')$ general
$\ell$-generated ideals in $IR'$ which are reductions of $IR'$. On
the other hand every dense open subset of
$\displaystyle\stackrel{t}{\times} {\mathbb G}(\ell, IR'
\otimes_{R'} K) =\displaystyle\stackrel{t}{\times} {\mathbb
G}(\ell, K^n)$ intersects with $\displaystyle\stackrel{t}{\times}
{\mathbb G}(\ell, I \otimes_R k)
=\displaystyle\stackrel{t}{\times} {\mathbb G}(\ell, k^n)$, since
$k$ is infinite. Thus ${\rm core}(IR') = {\mathfrak a}_1 \cap
\ldots \cap {\mathfrak a}_t$ where ${\mathfrak a}_i$ are minimal
reductions of $IR'$ extended from $R$-ideals ${\mathfrak b}_i$.
Now ${\mathfrak b}_i$ are minimal reductions of $I$, which gives
${\rm core}(I) \subset {\mathfrak b}_1 \cap \ldots \cap {\mathfrak
b}_t$. Hence $({\rm core}(I))R' \subset {\mathfrak b}_1 R' \cap
\ldots \cap {\mathfrak b}_t R' = {\rm core}(IR')$ and thus $({\rm
core}(I))R' = {\mathfrak b}_1 R' \cap \ldots \cap {\mathfrak b}_t
R' = {\rm core}(IR')$. In particular $({\rm core}(I))R' =
({\mathfrak b}_1 \cap \ldots \cap {\mathfrak b}_t)R'$, hence ${\rm
core}(I) = {\mathfrak b}_1 \cap \ldots \cap {\mathfrak b}_t$,
showing $\gamma(I) \leq \gamma(IR')$. The inequality $\gamma(I)
\geq \gamma(IR')$ is obvious since $({\rm core}(I))R' = {\rm
core}(IR')$. \QED \end{proof}

Lemmas~\ref{4.4} and~\ref{4.6} allow us to prove one of our main results:

\begin{theorem}\label{4.7}
Let $R$ be a local Cohen--Macaulay ring with infinite residue field
and $I$ an $R$-ideal of analytic spread $\ell$. Assume that $I$ is
$G_{\ell}$ and universally weakly $(\ell-1)$-residually $S_2$.
\begin{description}[(a)]
\item[$(${\it a}$)$]
Let ${\mathcal A}_1, \ldots, {\mathcal A}_t$
be universal $\ell$-generated ideals in $IS$. Then $({\rm core}(I))S
= {\rm core}(IS) = {\mathcal A}_1 \cap \ldots \cap {\mathcal
A}_t$ if and only if $t \geq \gamma(I) = \gamma(IS)$.

\item[$(${\it b}$)$]
Let ${\mathcal A}$ be a universal $\ell$-generated ideal in $IR(X)$.
Then ${\rm core}(I) = {\mathcal A} \cap R$.
\end{description}
\end{theorem}
\begin{proof} $(${\it a}$)$
If ${\rm core}(IS) = {\mathcal A}_1 \cap \ldots \cap
{\mathcal A}_t$ then $\gamma(IS) \leq t$ and so by Lemma~\ref{4.6},
$t \geq \gamma(I) = \gamma(IS)$.
Conversely, assume $t \geq \gamma(I)$. The equality $({\rm core}(I))S
= {\rm core}(IS)$ holds by Lemma~\ref{4.6}. As to the second equality,
Lemma~\ref{4.4} gives $e(IS/{\rm core}(IS)) =
e(IS/({\rm core}(I))S) \leq e(IS/{\mathcal A}_1 \cap \ldots \cap
{\mathcal A}_t)$. Now proceed as in the proof of Theorem~\ref{4.5}.

$(${\it b}$)$ In the setting of part $(${\it a}$)$, a subgroup of
${\rm Aut}_R(S)$ acts transitively on $\{ {\mathcal A}_1, \ldots,
{\mathcal A}_t \}$. Hence ${\mathcal A}_1 \cap R = \ldots =
{\mathcal A}_t \cap R$. Now $(${\it a}$)$ implies that ${\rm
core}(I)={\mathcal A}_1 \cap R$. We may assume that $R(X)
\hookrightarrow S$ is a flat local extension and ${\mathcal A}_1 =
{\mathcal A}S$, thus ${\mathcal A}_1 \cap R = {\mathcal A} \cap
R$. \QED
\end{proof}

Using Theorem~\ref{4.7}.{\it b} we can in turn generalize
Lemma~\ref{4.6} to the case when the map is not necessarily local:

\begin{theorem}\label{4.8}
Let $R \longrightarrow R'$ be a flat map of local Cohen--Macaulay rings
with infinite residue fields. Let $I$ be an $R$-ideal of analytic spread
$\ell$ such that $I$ and $IR'$ are $G_{\ell}$ and universally weakly
$(\ell-1)$-residually $S_2$. Then ${\rm core}(IR')=({\rm core}(I))R'$.
\end{theorem}
\begin{proof}
The containment ${\rm core}(IR') \subset ({\rm core}(I))R'$ is
obvious, because $\gamma(I) < \infty$ by Theorem~\ref{3.1}. If
$\ell(IR') < \ell$ then ${\rm core}(IR') = IR' \supset ({\rm
core}(I))R'$ by Lemma~\ref{2.2}.{\it g}. Thus we may assume that
$\ell(IR') = \ell$. Let ${\mathcal A}$ be a universal
$\ell$-generated ideal in $IR(X)$, $S=R(X)$, and $S' = R'(X)$.
Notice that ${\mathcal A}S'$ is a universal $\ell$-generated ideal
in $IS'$. Thus by Theorem~\ref{4.7}.{\it b}, ${\rm core}(IR') =
{\mathcal A}S' \cap R'$. Therefore $({\rm core}(IR')) \cap R =
{\mathcal A}S' \cap R' \cap R = {\mathcal A}S' \cap R = {\mathcal
A}S' \cap S \cap R \supset {\mathcal A} \cap R = {\rm core}(I)$,
where the last equality follows again from Theorem~\ref{4.7}.{\it
b}. Hence ${\rm core}(IR') \supset ({\rm core}(I))R'$. \QED
\end{proof}

\begin{corollary}\label{4.9}
Let $R$ be a local Cohen--Macaulay ring with infinite residue field and $I$
an $R$-ideal of analytic spread $\ell$. Assume that $I$ is $G_{\ell}$ and
universally weakly $(\ell-1)$-residually $S_2$. Then
\[
\gamma(I) = \max (\{ \gamma(I_{\mathfrak p}) \, | \, {\mathfrak p}
\in {\rm Min}({\rm Fitt}_{\ell}(I)) \} \cup \{ 1 \}).
\]
\end{corollary}
\begin{proof}
According to Lemma~\ref{2.2}.{\it a}, every localization of $I$ is
universally weakly $(\ell-1)$-residually $S_2$.
Write $t = \max (\{ \gamma(I_{\mathfrak p}) \, | \, {\mathfrak p} \in
{\rm Min}({\rm Fitt}_{\ell}(I)) \} \cup \{ 1 \})$, and notice that
$t$ and $\gamma(I)$ are finite by Theorem~\ref{3.1}.
For ${\mathfrak p} \in {\rm Spec}(R)$,
$({\rm core}(I))_{\mathfrak p} = {\rm core}(I_{\mathfrak p})$ by
Theorem~\ref{4.8} and hence $\gamma(I) \geq \gamma(I_{\mathfrak p})$.
Thus $\gamma(I) \geq t$. To prove the reverse inequality, let
${\mathcal A}_1, \ldots, {\mathcal A}_t$ be universal $\ell$-generated
ideals in $IS$ and recall that $\gamma(IS)$ is finite by
Theorem~\ref{4.7}.{\it a}. If ${\mathfrak p}
\in {\rm Min}({\rm Fitt}_{\ell}(I))$ then $\ell(I_{\mathfrak p})=\ell$
by Lemma~\ref{2.2}.{\it g}, and $({\mathcal A}_1)_{{\mathfrak p}S},
\ldots, ({\mathcal A}_t)_{{\mathfrak p}S}$ are universal $\ell$-generated
ideals in $I_{\mathfrak p} S_{{\mathfrak p}S}$. Thus by
Theorem~\ref{4.7}.{\it a}, ${\rm core}(I_{\mathfrak p} S_{{\mathfrak p}S}) =
({\mathcal A}_1)_{{\mathfrak p}S} \cap \ldots \cap ({\mathcal
A}_t)_{{\mathfrak p}S}$, and hence by Theorem~\ref{4.8},
$({\rm core}(IS))_{{\mathfrak p}S} = ({\mathcal A}_1 \cap \ldots \cap
{\mathcal A}_t)_{{\mathfrak p}S}$ for every ${\mathfrak p} \in
{\rm Min}({\rm Fitt}_{\ell}(I))$. Since furthermore
${\rm core}(IS) \subset {\mathcal A}_1 \cap \ldots \cap {\mathcal A}_t
\subset IS$ and ${\rm Ass}_S(IS/{\rm core}(IS)) \subset
{\rm Min}({\rm Fitt}_{\ell}(IS))$ by Lemma~\ref{2.2}.{\it e,{\,}f},
we conclude that
${\rm core}(IS) = {\mathcal A}_1 \cap \ldots \cap {\mathcal A}_t$. Now
$\gamma(I) \leq t$ by Theorem~\ref{4.7}.{\it a}. \QED
\end{proof}

\begin{remark}\label{4.10}
Except for the second assertion in Lemma~\ref{4.2}, the results of
this section remain true for ideals $I$ and $IR'$ primary to the
maximal ideals, if $R$ and $R'$ are merely Buchsbaum instead of
Cohen--Macaulay. Notice that the first part of Lemma~\ref{4.2}
still holds in this case, because ${\rm length}(I/J)$ is
independent of the minimal reduction of $I$ by definition of the
Buchsbaum property.
\end{remark}

However, even for an ${\mathfrak m}$-primary ideal $I$ the inclusion
${\mathcal A} \cap R \subset {\rm core}(I)$ of Theorem~\ref{4.7}.{\it b}
may fail to hold if $(R, {\mathfrak m})$ is not Buchsbaum, as can be
seen from \cite[p. 246]{RS}. We observe a similar failure for
Cohen--Macaulay rings provided the ideal $I$ does not satisfy the
assumptions of Theorem~\ref{4.7}.

\smallskip

Let $R$ be a Noetherian local ring with infinite residue field and $I$ an
$R$-ideal. If $J$ is a reduction of $I$ we denote by $r_J(I)$ the least
integer $r \geq 0$ with $I^{r+1}=JI^r$. Recall that the {\it reduction number}
of $I$ is defined as $r(I)= \min \{ r_J(I) \, | \, J \in {\mathcal M}(I) \}$.

\begin{example}\label{4.11}
Let $R = k[U, V, W]_{(U, V, W)}/(U^2+V^2, VW)$, where $k$ is an infinite
field and $U, V, W$ are variables. Denote the images of $U, V, W$ in $R$
by $u, v, w$. Consider the $R$-ideal $I=(u, v)$, and let ${\mathcal A} =
(Xu+Yv) \subset S = R(X, Y)$ be a universal one-generated ideal in $IS$.

Notice that $R$ is Gorenstein, $\ell = \ell(I) = 1$, and $I$ does
not satisfy $G_1$, but is universally $0$-residually $S_2$. In
this case ${\rm core}(I)$ is the intersection of finitely many
minimal reductions of $I$, ${\rm core}(IR') = ({\rm core}(I))R'$
for every flat map $R \longrightarrow R'$ to a local
Cohen--Macaulay ring $R'$, but ${\rm core}(I)$ is not an
intersection of general one-generated ideals in $I$ which are
reductions of $I$, and ${\rm core}(I) \subsetneq {\mathcal A} \cap
R$.

Indeed, $(u)$ and $(v)$ are minimal reductions of $I$, hence ${\rm
core}(I) \subset (u) \cap (v) = I^2$. On the other hand the
special fiber ring ${\rm gr}_I(R) \otimes_R k$ is defined by a
single quadric; hence $r_J(I) = 1$ for every minimal reduction $J$
of $I$. Thus $I^2 \subset {\rm core}(I)$. Therefore
\[
{\rm core}(I) = I^2 = (u) \cap (v).
\]
If the map $R \longrightarrow R'$ is local then the same argument gives
${\rm core}(IR') = (IR')^2$. Otherwise either ${\rm core}(IR') = 0 =
(IR')^2$ or ${\rm core}(IR') = R' = (IR')^2$.
Hence in any case ${\rm core}(IR') = ({\rm core}(I))R'$.

A general one-generated ideal $(\lambda u + \mu v)$ in $I$
contains $uw = (\lambda u + \mu v) \lambda^{-1} w$, whereas $uw
\not\in I^2$. Thus ${\rm core}(I)$ cannot be the intersection of
general one-generated ideals in $I$. Likewise $I^2 + (uw) \subset
{\mathcal A}$, hence ${\rm core}(I) \subsetneq {\mathcal A} \cap
R$.
\end{example}

\section{Computational remarks}

Individual minimal reductions of homogeneous ideals tend to be
inhomogeneous -- for instance, the monomial ideal $I=(U^2, UV, V^3)
\subset k[U, V]_{(U, V)}$ has no minimal reduction generated by
homogeneous polynomials in $U$ and $V$. Nevertheless the core of
this ideal is monomial due to the following general fact:

\begin{remark}
Let $R$ be a Noetherian local ring with infinite residue field.
Assume that $R=R'_{{\mathfrak m}'}$ for $R'$ an ${\mathbb
N}_0^n$-graded ring over a local ring and ${\mathfrak m}'$ its
homogeneous maximal ideal. Let $I$ be an $R$-ideal.
If $I$ is generated by homogeneous elements of $R'$ then so is
${\rm core}(I)$.
\end{remark}
\begin{proof}
$($This proof was suggested to us by D. Eisenbud.$)$ Let $U$ be
the group of units of $[R']_{(0, \ldots, 0)}$ and $G$ the direct
product $U^n$. If $\alpha = (\alpha_1, \ldots, \alpha_n) \in G$
and $x \in [R']_{(i_1, \ldots, i_n)}$ we define $\alpha x$ to be
$\alpha_1^{i_1} \cdots \alpha_n^{i_n} x \in [R']_{(i_1, \ldots,
i_n)}$. This induces an action of $G$ on the ring $R$. As is well
known, an $R$-ideal is $G$-stable if and only if it is extended
from a homogeneous $R'$-ideal. To finish the proof, notice that
$G$ acts on the set ${\mathcal M}(I)$, which implies the
$G$-stability of ${\rm core}(I)$. \QED \end{proof}

The next remark gives a fairly efficient probabilistic algorithm for
computing the core. In light of Theorem~\ref{4.5}
it suffices to bound $\gamma(I)$:

\begin{remark}\label{5.2}
Let $R$ be a local Cohen--Macaulay ring with infinite residue field
and $I$ an $R$-ideal of analytic spread $\ell$.
Assume that $I$ is $G_{\ell}$ and weakly $(\ell-1)$-residually $S_2$.
If ${\mathfrak a}_1 \cap \ldots \cap {\mathfrak a}_t = {\mathfrak a}_1
\cap \ldots \cap {\mathfrak a}_{t+1}$ for ${\mathfrak a}_1,
\ldots, {\mathfrak a}_{t+1}$ general $\ell$-generated ideals in $I$,
then $\gamma(I) \leq t$.
\end{remark}
\begin{proof}
Let $k$ be the residue field of $R$. By Theorem~\ref{3.1} and
Theorem~\ref{4.5} there exists an integer $s > t$ so that ${\rm
core}(I) = {\mathfrak a}_1 \cap \ldots \cap {\mathfrak a}_s$ with
${\mathfrak a}_1, \ldots, {\mathfrak a}_s$ general
$\ell$-generated ideals in $I$. After passing to a smaller dense
open subset of $\displaystyle\stackrel{s}{\times} {\mathbb
G}(\ell, I \otimes_R k)$, we deduce from our assumption that
${\mathfrak a}_1 \cap \ldots \cap {\mathfrak a}_t = {\mathfrak
a}_1 \cap \ldots \cap {\mathfrak a}_t \cap {\mathfrak a}_i$ for
every $t+1 \leq i \leq s$. Thus ${\mathfrak a}_1 \cap \ldots \cap
{\mathfrak a}_t = {\mathfrak a}_1 \cap \ldots \cap {\mathfrak a}_s
= {\rm core}(I)$. \QED \end{proof}

As before, it suffices to assume in Remark~\ref{5.2} that $R$ is
Buchsbaum if $I$ is primary to the maximal ideal.

\begin{lemma}\label{5.3}
Let $R$ be a local Cohen--Macaulay ring, $s$ an integer,
$I$ an $R$-ideal satisfying $G_s$, ${\mathcal B}$ a generic
$s$-generated ideal in $IT$, $K={\mathcal B} :_T IT$, and
${\mathfrak q} \in V(K)$ with ${\rm dim}\, T_{\mathfrak q} \leq s$.

\begin{description}[(a)]
\item[$(${\it a}$)$]
If $I \subset {\mathfrak q}$ then ${\mathfrak q}$ is extended from
a minimal prime of \/ ${\rm Fitt}_s(I)$.

\item[$(${\it b}$)$]
If $I \not\subset {\mathfrak q}$ then ${\mathfrak q} \cap R$
is a minimal prime of $R$.
\end{description}
\end{lemma}
\begin{proof}
We write ${\mathfrak p} ={\mathfrak q} \cap R$ and replace $R$ by
$R_{\mathfrak p}$.

$(${\it a}$)$ Notice that ${\rm dim}\, R \leq {\rm dim}\,
T_{\mathfrak q} \leq s \leq {\rm ht}\, {\rm Fitt}_s(I)$, where the
last inequality is a consequence of the $G_s$ assumption. Thus it
suffices to prove that ${\rm Fitt}_s(I) \not= R$, since then ${\rm
dim}\, T_{\mathfrak q} = {\rm ht}\, {\rm Fitt}_s(I)$. Suppose
${\rm Fitt}_s(I) = R$. In this case $I$ satisfies $G_{s+1}$, and
hence ${\rm ht} (IT+K) \geq s+1$ by \cite[3.2]{HU}. But this is
impossible because $IT+K \subset {\mathfrak q}$ and ${\rm dim}\,
T_{\mathfrak q} \leq s$.

$(${\it b}$)$ One has ${\mathcal B} \subset {\mathfrak q}$ since
${\mathfrak q} \in V(K)$, and $I=R$ since $I \not\subset
{\mathfrak q}$. Now after adjoining variables to $T$ and applying
an $R$-automorphism we may suppose that ${\mathcal B}$ is defined
using $1$ as a generator of $I$. Hence ${\mathcal B}$ is generated
by $s$ variables $X_1, \ldots, X_s$ of $T$, and thus $({\mathfrak
p}, X_1, \ldots, X_s) \subset {\mathfrak q}$. As ${\rm dim}\,
T_{\mathfrak q} \leq s$, ${\mathfrak p}$ must be a minimal prime
of $R$. \QED \end{proof}

\begin{proposition}\label{5.4}
Let $R$ be a local Cohen--Macaulay ring with infinite
residue field and $I$ an $R$-ideal of analytic spread $\ell$
and reduction number $r$.
Let $f \in I$ and $h \in \sqrt{{\rm Fitt}_{\ell}(I)}$ be non zerodivisors
on $R$, and ${\mathcal B}$ a generic $\ell$-generated ideal in $IT$.
Assume that $I$ is $G_{\ell}$ and universally weakly $(\ell-1)$-residually
$S_2$, and $IT$ is weakly $(\ell-1)$-residually $S_2$.
Then
\[
{\rm core}(I) = [{\mathcal B} :_T ({\mathcal B} :_T h^{\infty}I)]_0
= [{\mathcal B} :_T ({\mathcal B} :_T f^{r+1})]_0.
\]
\end{proposition}
\begin{proof}
Write ${\mathfrak m}$ for the maximal ideal of $R$, $S =
T_{{\mathfrak m}T}$ and ${\mathcal A} = {\mathcal B}S$. By
Theorem~\ref{4.7}.{\it b}, ${\rm core}(I) = {\mathcal A} \cap R$.
Now ${\mathcal A} \cap R = {\mathcal A} \cap T \cap R =[{\mathcal
B}_{{\mathfrak m}T} \cap T]_0$. Write $H = {\mathcal B} \colon
({\mathcal B} \colon h^{\infty}I)$ and $F = {\mathcal B} \colon
({\mathcal B} \colon f^{r+1})$. It remains to prove that
${\mathcal B}_{{\mathfrak m}T} \cap T = H = F$. This will follow
once we have shown that $H_{{\mathfrak m}T} = F_{{\mathfrak m}T} =
{\mathcal B}_{{\mathfrak m}T}$ and that every associated prime of
$H$ or $F$ is contained in ${\mathfrak m}T$.

First, notice that ${\mathcal B}_{{\mathfrak m}T}$ is an $\ell$-generated
reduction of $IT_{{\mathfrak m}T}$. Hence by Lemma~\ref{2.2}.{\it e},
$({\mathcal B} \colon IT)_{{\mathfrak m}T}$ contains some power of
${\rm Fitt}_{\ell}(I)$ and hence of $h$. This gives
$({\mathcal B} \colon h^{\infty}I)_{{\mathfrak m}T} =
(({\mathcal B} \colon IT) \colon h^{\infty})_{{\mathfrak m}T}
= T_{{\mathfrak m}T}$. Therefore
$H_{{\mathfrak m}T} = {\mathcal B}_{{\mathfrak m}T}$.
Since ${\mathcal B}_{{\mathfrak m}T}$ is a universal $\ell$-generated
ideal in $IT_{{\mathfrak m}T}$ we have $I^{r+1} \subset {\mathcal
B}_{{\mathfrak m}T}$ by \cite[3.4]{SUV}, and hence $f^{r+1} \in
{\mathcal B}_{{\mathfrak m}T}$. This gives
$({\mathcal B} : f^{r+1})_{{\mathfrak m}T} = T_{{\mathfrak m}T}$,
thus $F_{{\mathfrak m}T} = {\mathcal B}_{{\mathfrak m}T}$.

Finally let ${\mathfrak q}$ be an associated prime of $H$ or $F$. Notice that
${\mathfrak q}$ is also an associated prime of ${\mathcal B}$.
Since ${\rm ht}\, {\mathcal B} \colon IT \geq \ell$ by \cite[the proof
of~3.2]{HU}, Lemma~\ref{2.2}.{\it c} then gives
${\rm dim} \, T_{\mathfrak q} \leq \ell$.
We claim that ${\mathfrak q} \subset {\mathfrak m}T$.
We may assume that ${\mathcal B} \colon IT \subset {\mathfrak q}$ since
otherwise $IT_{\mathfrak q} = {\mathcal B}_{\mathfrak q}$, thus
${\mathfrak q}$ is an associated prime of $IT$ and hence contained in
${\mathfrak m}T$. Now if $I \subset {\mathfrak q}$ then ${\mathfrak q} \subset
{\mathfrak m}T$ by Lemma~\ref{5.3}.{\it a}. If on the other hand
$I \not\subset {\mathfrak q}$ then part $(${\it b}$)$ of the same lemma
gives $h \not\in {\mathfrak q}$ and $f \not\in {\mathfrak q}$.
Therefore $({\mathcal B} : h^{\infty}I)_{\mathfrak q} =
{\mathcal B}_{\mathfrak q}$ and $({\mathcal B} : f^{r+1})_{\mathfrak q} =
{\mathcal B}_{\mathfrak q}$. Thus $H_{\mathfrak q} = F_{\mathfrak q} =
T_{\mathfrak q}$, which is impossible. \QED
\end{proof}

\begin{remark}\label{5.5}
Let $(R, {\mathfrak m})$ be a local Cohen--Macaulay ring of
dimension $d$ with infinite residue field, $I$ an ${\mathfrak
m}$-primary ideal of multiplicity $e$, $f$ an $R$-regular element,
and ${\mathcal B}$ a generic $d$-generated ideal in $IT$. Then
\[
{\rm core}(I) = [{\mathcal B} :_T ({\mathcal B} :_T f^e)]_0.
\]
\end{remark}
\begin{proof}
The assertion follows from the proof of Proposition~\ref{5.4} and the
fact that $\lambda((T/{\mathcal B})_{{\mathfrak m}T})=e$. \QED
\end{proof}

The equality of Proposition~\ref{5.4} gives a method for computing
the core of a broad class of ideals generated by homogeneous
polynomials not necessarily of the same degree: by giving suitable
degrees to the variables $X_{jl}$ of $T$, the ideal ${\mathcal B}$
becomes homogeneous and the computation stays in the graded
category. As an illustration, taking $I=(U^3, UV^3, V^4) \subset
k[ U, V ]_{(U, V)}$, we obtain ${\rm core}(I)=(U^2, UV, V^2)I$ and
taking $I = (U^3, UV^2W^2, V^3W^3) \subset k[U,V,W]_{(U,V,W)}$, we
obtain ${\rm core}(I) = (U^2, UVW, V^2W^2)I$. The outcome of
neither computation could have been predicted by the results of
\cite{CPU} or \cite{HS}.


\bigskip

\begin{thebibliography}{99}

\bibitem{AH}{L. Avramov and J. Herzog, The Koszul algebra of a
codimension $2$ embedding, Math. Z. {\bf 175} (1980), 249-260}

\bibitem{CEU}{M. Chardin, D. Eisenbud and B. Ulrich, Hilbert functions,
residual intersections, and residually $S_2$-ideals, Compositio
Math. {\bf 125} (2001), 193-219}

\bibitem{CPU}{A. Corso, C. Polini and B. Ulrich, Core and residual
intersections of ideals, preprint 2001}

\bibitem{EHV}{D. Eisenbud, C. Huneke and W.V. Vasconcelos, Direct methods
for primary decomposition, Invent. Math. {\bf 110} (1992),
207-235}

\bibitem{HVV}{J. Herzog, W.V. Vasconcelos and R.H. Villarreal,
Ideals with sliding depth, Nagoya Math. J. {\bf 99} (1985),
159-172}

\bibitem{Huneke}{C. Huneke, On the symmetric and Rees algebra of an ideal
generated by a $d$-sequence, J. Algebra {\bf 62} (1980),
268-275.}

\bibitem{H1}{C. Huneke, Linkage and Koszul homology of ideals,
Amer. J. Math. {\bf 104} (1982), 1043-1062.}

\bibitem{HS}{C. Huneke and I. Swanson, Cores of ideals in $2$-dimensional
regular local rings, Michigan Math. J. {\bf 42} (1995), 193-208.}

\bibitem{HU}{C. Huneke and B. Ulrich, Residual intersections, J. reine
angew. Math. {\bf 390} (1988), 1-20.}

\bibitem{L}{J. Lipman, Adjoints of ideals in regular local rings,
Math. Research Letters {\bf 1} (1994), 1-17.}

\bibitem{LS}{J. Lipman and A. Sathaye, Jacobian ideals and a theorem
of Brian\c{c}on--Skoda, Michigan Math. J. {\bf 28} (1981), 199-222.}

\bibitem{Nagata}{M. Nagata, {\it Local Rings}, Krieger, Huntington,
1975.}

\bibitem{RS}{D. Rees and J.D. Sally, General elements and joint
reductions, Michigan Math. J. {\bf 35} (1988), 241-254.}

\bibitem{SUV}{A. Simis, B. Ulrich, and W.V. Vasconcelos,
Cohen--Macaulay Rees algebras and degrees of polynomial relations,
Math. Ann. {\bf 301} (1995), 421-444}

\bibitem{AN}{B. Ulrich, Artin-Nagata properties and reductions of ideals,
in Commutative Algebra: Syzygies, Multiplicities, and Birational
Algebra, W. Heinzer, C. Huneke, J. Sally (eds.), Contemp. Math.
{\bf 159}, Amer. Math. Soc., Providence, 1994, 373-400}

\bibitem{Valla}{G. Valla, On the symmetric and Rees algebras of an ideal,
Manuscripta Math. {\bf 30} (1980), 239-255}

\bibitem{Vas}{W.V. Vasconcelos, Reduction numbers of ideals,
J. Algebra {\bf 216} (1999), 652-664}

\end{thebibliography}
\end{document}